\documentclass[10pt]{amsart}

\usepackage{amsmath,amssymb,amsfonts,epsfig}
\usepackage{hhline}

\newtheorem{theorem}{Theorem}[section]

\theoremstyle{definition}
\newtheorem{definition}[theorem]{Definition}

\theoremstyle{remark}
\newtheorem{remark}[theorem]{Remark}
\newtheorem{example}[theorem]{Example}

\numberwithin{equation}{section}

\def\wo{\overline}

\def\R{{\mathbb R}}
\def\Z{{\mathbb Z}}
\def\a{\alpha}
\def\b{\beta}
\def\d{\delta}
\def\Ga{\Gamma}
\newcommand{\rb}{\raisebox}
\newcommand{\ig}{\includegraphics}
\newcommand\risS[6]{\rb{#1pt}[#5pt][#6pt]{\begin{picture}(#4,15)(0,0)
  \put(0,0){\ig[width=#4pt]{#2.eps}} #3
     \end{picture}}}

\input xy
\xyoption{all}

\begin{document}

\title[Topological Tutte Polynomial]
  {Topological Tutte Polynomial}
\author{SERGEI CHMUTOV}
\date{}

\address{Ohio State University - Mansfield, 
1760 University Drive, 
Mansfield, OH 44906, USA \linebreak 
{\tt chmutov@math.ohio-state.edu}}

\keywords{Graphs on surfaces, ribbon graphs, Bollob\'as-Riordan polynomial, Tutte polynomial, Las Vergnas polynomial, Krushkal polynomial, duality, virtual links, quasi-trees.}

\begin{abstract}
This is a survey of recent works on topological extensions of the Tutte polynomial. 
\end{abstract}

\maketitle
{\center\small Mathematics Subject Classifications: 05C31, 05C22, 05C10, 57M15, 57M25, 57M30}

\section{Introduction}  
We survey recent works on topological extensions of the Tutte polynomial. 

The extensions deal with graphs on surfaces. Graphs embedded into surfaces are the main objects of the topological graph theory. There are several books devoted to this subject \cite{MR3086663,MR1855951,MR2036721,MR1844449}.
Cellular embeddings of graphs can be described as {\em ribbon graphs}. Oriented ribbon graphs appear under different names such as {\it rotation systems} \cite{MR1844449}, {\it maps}, {\it fat graphs}, {\it cyclic graphs}, {\it dessins d'enfants} \cite{MR2036721}. 
Since the pioneering paper of L.~Heffter in 1891 \cite{He} they occur in various parts of mathematics ranging from graph theory, combinatorics, and topology to representation theory, Galois theory, algebraic geometry, and quantum field theory \cite{MR2962302,MR0830043,MR1341841,MR2036721,MR0919235,MR1036112}. 
For example, ribbon graphs are used to enumerate cells in the cell decomposition of the moduli spaces of complex algebraic curves
\cite{MR0830043,MR0919235,MR2036721}. The absolute Galois group 
$\mathrm{Aut}(\wo{\mathbb Q}/{\mathbb Q})$ faithfully acts on the set of ribbon graphs (see \cite{MR2036721} 
and references therein). Ribbon graphs are the main combinatorial objects of the Vassiliev knot invariant theory 
\cite{MR2962302}. 
They are very useful for Hamiltonicity of the Cayley graphs \cite{MR2341831,MR1417567}.

M.~Las Vergnas \cite{MR0597150} found a generalization of the Tutte polynomial to cellularly embedded graphs as an application of his matroid perspectives. Recently his polynomial was extended to not necessarily cellularly embedded graphs \cite{EMM1}.

In 2001 B.~Bollob\'as and O.~Riordan \cite{MR1851080,MR1906909}, motivated by some problems of knot theory, introduced a different generalization of the Tutte polynomial for ribbon graphs. For non-planar graphs, there is no duality relation for the Tutte polynomial but there is one for the Bollob\'as-Riordan polynomial. In \cite{MR1906909}, this was proved for a one-variable specialization. 
J.~Ellis-Monaghan and I.~Sarmiento \cite{EMS} extended it to a two-variable relation. 
Independently the duality formula was discovered in \cite{MR2368618}.

In knot theory the classical Thistlethwaite's theorem relates the Jones polynomial and the Tutte 
polynomial of a corresponding planar graph \cite{MR0935433,MR0899051}. 
Igor Pak suggested using the Bollob\'as-Riordan polynomial for Thistlethwaite type theorems. This idea was first realized in \cite{MR2343139} for a special class of (checkerboard colorable) virtual links. Then it was realized for classical links in \cite{MR2389605}, and for arbitrary virtual links in \cite{MR2460170}. 
Formally all three theorems from \cite{MR2343139,MR2389605,MR2460170} were different. They used different constructions of a ribbon graph from a link diagram and different substitutions in the Bollob\'as-Riordan polynomials of these graphs. An attempt to understand and unify these theorems led to the discovery of {\it partial duality} in \cite{MR2507944} (called there {\it generalized duality}) and to a proof of an invariance of a certain specializations of the 
Bollob\'as-Riordan polynomial under it.

A different generalization of Thistlethwaite theorem for virtual
links was found in \cite{MR3111544} based on the relative Tutte polynomial of plane graphs. The equivalence of this approach to that of \cite{MR2507944} was clarified in \cite{MR3390227}.

In 2011 V.~Krushkal \cite{MR2769192} found a four variable generalization of the Tutte polynomial. It can be reduced to the Bollob\'as-Riordan polynomial under a certain substitution. It turns out \cite{MR2998838} that it also can be reduced to the Las Vergnas polynomial under a different substitution, while the Las Vergnas and the Bollob\'as-Riordan polynomials are essentially independent of each other. Generalizing the classical spanning tree expansion of the Tutte polynomial, Clark Butler found \cite{But} an elegant quasi-tree expansion of the Krushkal polynomial.

The next diagram represents various relations between these polynomials.
$$
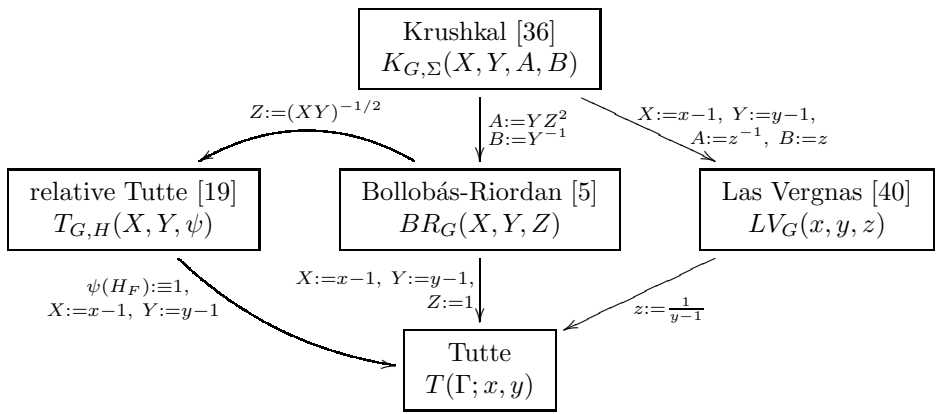
\xymatrix{
&\fbox{$\begin{array}{c}\mbox{Krushkal \cite{MR2769192}}\\
        K_{G,\Sigma}(X,Y,A,B)\end{array}$}
   \ar[d]^{\hspace*{-5pt}\begin{array}{l}\scriptstyle               
                                     A:=YZ^2\vspace{-5pt}\\
             \scriptstyle B:=Y^{-1}\end{array}}
   \ar[rd]^(.7){\hspace*{-20pt}\begin{array}{l}\scriptstyle
                                  X:=x-1,\ Y:=y-1,\vspace{-3pt}\\
             \scriptstyle\hspace*{20pt}A:=z^{-1},\ B:=z\end{array}} &\\
\fbox{$\begin{array}{c}\mbox{relative Tutte \cite{MR3111544}}\\
       T_{G,H}(X,Y,\psi)\end{array}$} 
 \ar@/_1.5pc/[rd]^{\hspace*{-105pt}\begin{array}{l}\scriptstyle
             \scriptstyle\hspace*{15pt}\psi (H_F):\equiv 1, 
                         \vspace{-3pt}\\
             \scriptstyle X:=x-1,\ Y:=y-1
             \end{array}} &
\fbox{$\begin{array}{c}\mbox{Bollob\'as-Riordan \cite{MR1906909}}\\
       BR_G(X,Y,Z)\end{array}$}
 \ar@/_2.5pc/[l]_{\rb{-18pt}{$\ \ \scriptstyle Z:=(XY)^{-1/2}$}} 
 \ar[d]_{\begin{array}{r}\scriptstyle
                                  X:=x-1,\ Y:=y-1,\vspace{-3pt}\\
             \scriptstyle Z:=1\end{array}\hspace*{-6pt}} & 
\fbox{$\begin{array}{c}\mbox{Las Vergnas \cite{MR0597150}}\\
       LV_G(x,y,z)\end{array}$} 
   \ar[ld]^{\scriptstyle z:=\frac1{y-1}} \\
&\fbox{$\begin{array}{c}\mbox{Tutte}\\
       T(\Ga;x,y)\end{array}$} &
}
$$

Both the relative Tutte polynomial of \cite{MR3111544} and the
Las Vergnas polynomial of \cite{MR0597150} may be formulated for matroids. 
Since the Bollob\'as-Riordan polynomial specializes to the relative Tutte polynomial one may expect that the relative Tutte and the Las Vergnas polynomials are also independent. This may signify the existence of a more general matroid polynomial which would be a matroidal counterpart of the Krushkal polynomial. 
Very recently I.~Moffatt and B.~Smith found such a polynomial.


The Tutte polynomial has also been extended to higher dimensional cell complexes. We do not pursue this connection here, but the interested reader may see  \cite{KR}. 
In \cite{MR3157807} this polynomial was related to the cellular spanning trees of 
\cite{MR2794024,MR0720308},
the arithmeric matroids of \cite{MR2989987} and the simplicial chromatic polynomial 
\cite{MR1215962,MR1304401}. 

\section{Ribbon graphs.}

The topological objects we are dealing with are graphs embedded into a surface. The embedding is assumed to be {\it cellular}, that means each connected component of the complement of the graph, a {\it face}, is homeomorphic to a disc.
\index{cellular embedding}\index{embedding!cellular} 

Considering a small regular neighborhood of the graph on a surface we come to the equivalent notion of {\it ribbon graphs}.

\begin{definition}[\cite{MR1851080,MR1906909}]
A {\it ribbon graph}\index{ribbon graph} 
is an abstract (not necessarily orientable) compact surface with 
boundary, decomposed into a number of closed topological discs of two types, 
{\it vertex-discs} and {\it edge-ribbons}, satisfying the following natural conditions: 
the discs of the same type are pairwise disjoint; the vertex-discs and the edge-ribbons intersect in disjoint line segments, each such line segment lies on the boundary of precisely one vertex and precisely one edge, and every edge contains exactly two such line segments.
\end{definition}
We consider ribbon graphs up to a homeomorphism of the corresponding surfaces preserving the decomposition into vertex-discs and edge-ribbons.
It is important to note that a ribbon graph is an abstract two-dimensional surface with boundary; its embedding into the 3-space shown in pictures is irrelevant.

Here are three examples.
$$\mbox{(a)}
\risS{-35}{rg-ex1}{}{40}{0}{0}\quad
  \mbox{(b)}
\risS{-30}{rg-ex2}{}{40}{0}{0}\quad
  \mbox{(c)}
\risS{-30}{rg-ex3}{}{180}{0}{45}
$$

Gluing a disc to each boundary component of a ribbon graph $G$ we get a closed surface without boundary. Placing a vertex at the center of each vertex-disc and connecting the vertices by edges along the edge-ribbons we get a graph, the {\it core} graph of $G$, cellularly embedded into the surface. Thus ribbon graphs and cellularly embedded graphs are the same objects.

\bigskip
In this paper we will use the following notation for the major parameters of a ribbon graph.
\begin{definition}
\begin{itemize}
\item[$\bullet$] $v(G) := |V(G)|$ denotes the number of vertices of a ribbon graph $G$;
\item[$\bullet$] $e(G) := |E(G)|$ denotes the number of edges of $G$;
\item[$\bullet$] $k(G)$ denotes the number of connected components of $G$;
\item[$\bullet$] $r(G):=v(G)-k(G)$ denotes the {\it rank} of $G$;
\item[$\bullet$] $n(G):=e(G)-r(G)$ denotes the {\it nullity} of $G$;
\item[$\bullet$] $bc(G)$ denotes the number of connected components of the boundary of the surface of $G$.
\end{itemize} 

Only the last parameter $bc(G)$ is a topological one. All the others parameters are standard in graph theory.

A {\it spanning subgraph} $F$ of a ribbon graph $G$ is 
a ribbon graph consisting of all the vertices of $G$ and a subset of the edges of $G$.
\end{definition}

\section{The Bollob\'as--Riordan polynomial.}\label{s:BRpol} 

The Bollob\'as--Riordan polynomial, originally defined in \cite{MR1851080,MR1906909}, was generalized to a multivariable polynomial of weighted ribbon graphs in
\cite{MR2368618, MR2552629}. 
We will use a slightly more general doubly weighted 
Bollob\'as--Riordan polynomial of a ribbon graph $G$ with weights $(x_e,y_e)$ of an edge $e\in G$.
The doubly weighted Tutte polynomial is useful in connection with the knot theory 
\cite{MR0955462}. Double weights also allow to specialize the Bollob\'as--Riordan polynomial to different signed versions of the Tutte polynomial, see the Section \ref{ssec:signed-BR}.

\bigskip
\begin{definition}\index{Bollob\'as--Riordan polynomial}
$$BR_G(X,Y,Z):=\sum_{F\subseteq G} (\prod_{e\in F}x_e) 
(\prod_{e\not\in F}y_e) X^{r(G)-r(F)}Y^{n(F)}Z^{k(F)-bc(F)+n(F)},
$$
where the sum runs over all spanning subgraphs $F$. 
\end{definition}

The original Bollob\'as--Riordan polynomial of \cite{MR1906909} used an
additional variable $w$ responsible for orientability of the ribbon subgraph $F$. We set $w=1$ and replace their $x-1$ by our $X$, but we introduce weights of edges.

Note that the exponent $k(F)-bc(F)+n(F)$ of the variable $Z$ is equal to
$2k(F)-\chi(\widetilde{F})$, where $\chi(\widetilde{F})$ is the Euler characteristic of the surface $\widetilde{F}$ obtained by gluing a disc to each boundary component of $F$.
For orientable $F$, it is twice the genus of $F$. For non-orientable $F$, according to the  classical classification theorem for surfaces, $\widetilde{F}$ can be represented by a connected sum of several real projective planes; the number of them is equal to the exponent of $Z$. In particular, for a planar ribbon graph $G$ (i.e. when the surface $G$ has genus zero) the 
Bollob\'as-Riordan polynomial $BR_G$ does not contain $Z$. 
In this case, and if all weights are equal to 1, it is equal to the classical Tutte polynomial 
$T(\Ga;x,y)$ of the core graph $\Ga$ of $G$ with shifted variables $X=x-1$, $Y=y-1$:
$$BR_G(x-1,y-1,Z) = T(\Ga;x,y)\ .$$
Similarly, a specialization $Z=1$ (and $x_e=1, y_e=1$ for all edges $e$) of the Bollob\'as-Riordan polynomial
of an arbitrary ribbon graph $G$ gives the Tutte polynomial of the core graph:
$$BR_G(x-1,y-1,1) = T(\Ga;x,y)\ .$$
So one may think of the Bollob\'as-Riordan polynomial as a generalization of the Tutte polynomial to graphs cellularly embedded into a surface and capturing topological information of the embedding.

\begin{example}
Let the following ribbon graph $G$ have all weights equal to 1.
$$
\begin{array}{c||c|c|c|c} \label{br-table}
\risS{8}{rg-ex32}{\put(7,-5){$G$}}{50}{30}{5}
& \quad\risS{8}{rgBBB}{}{45}{0}{0}
& \quad\risS{8}{rgBBA}{}{45}{0}{0} 
& \quad\risS{13}{rgBAB}{}{45}{0}{0}  
& \quad\risS{13}{rgBAA}{}{45}{0}{0}\\ \hline
 (k,r,n,bc)\atop {\rm term\ of\ } BR_G & (1,1,1,2)\atop Y &
  (1,1,0,1)\atop 1 &
  (1,1,0,1)\atop 1 & {(2,0,0,2)\atop X} \rb{-10pt}{\makebox(0,25){}}
\\ \hline\hline
& \risS{8}{rg-ex32}{}{48}{40}{5}
& \risS{8}{rgABA}{}{48}{0}{0} 
& \risS{13}{rgAAB}{}{48}{0}{0}
& \risS{13}{rgAAA}{}{48}{0}{0}\\ \cline{2-5}
& (1,1,2,1)\atop Y^2Z^2 & (1,1,1,1)\atop YZ & (1,1,1,1)\atop YZ &
 {(2,0,1,2)\atop XYZ} \rb{-10pt}{\makebox(0,25){}}
\end{array}
$$
$$BR_G(X,Y,Z) = Y+2+X+Y^2Z^2+2YZ+XYZ$$
\end{example}

\subsection{Properties.}
One can easily prove
\begin{theorem}
$$BR_G=\left\{\begin{array}{ll}
x_e\cdot BR_{G/e} + y_e\cdot BR_{G-e} & 
         \mbox{if $e$ is ordinary, that is neither}\\
         &\mbox{a bridge nor a loop,} \\
x_e\cdot BR_{G/e} + y_eX\cdot BR_{G-e} & \mbox{if $e$ is a bridge.} \\
x_eY\cdot BR_{G/e}+y_e\cdot BR_{G-e} & \mbox{if $e$ is a trivial orientable loop,}\\
x_eYZ\cdot BR_{G/e}+y_e\cdot BR_{G-e} & \mbox{if $e$ is a non-orientable loop.}
\end{array}\right.
$$
$$
BR_{G_1\sqcup G_2} = BR_{G_1}\cdot BR_{G_2}, 
\mbox{where $\sqcup$ is disjoint union.}
$$
\end{theorem}
Here a {\it trivial loop} means a loop such that the deleting it and cutting its vertex-disc along a chord connecting the endpoints of the loop separates the the ribbon graph.

\subsection{A signed version of the Bollob\'as-Riordan polynomial.}\label{ssec:signed-BR}
Signed graphs are weighted graphs with single weights of $\pm1$ on their edges. They are very important in particular in knot theory. There is a vast literature on signed graphs 
\cite{MR1744869}. 

In \cite{MR1031266} 
L.~Kauffman introduced a generalization of the Tutte polynomial to signed graphs. The signed version of the Bollob\'as-Riordan polynomial was introduced in \cite{MR2343139} 
and used in \cite{MR2507944} 
(a version of it was also used in \cite{MR2382735}). 
This signed version can be obtained by changing the sign weights to double weights and using the Bollob\'as-Riordan polynomial above. We set the weights of positive edges to 1, $x_+:=y_+:=1$, and the weights of negative edges to $x_-:=(X/Y)^{1/2}$, $y_-:=(Y/X)^{1/2}$. Thus the signed version becomes
$$BR_G(x,y,z)\ =\ \sum_{F \in G} X^{r(G)-r(F)+s(F)} Y^{n(F)-s(F)}
   Z^{k(F)-bc(F)+n(F)}\ , 
   \index{Bollob\'as--Riordan polynomial!signed version}
$$
where
$s(F) := \frac{e_{-}(F)-e_{-}(\wo F)}{2}$,
$e_{-}(F)$ denotes the number of negative edges in $F$, and 
$\wo F = G-F$ is the complement to $F$ in $G$, i.e.
the spanning subgraph of $G$ with exactly those (signed) edges of $G$ that do not belong to $F$.

As above, specializing to $Z=1$ gives Kauffman's Tutte polynomial of signed graphs \cite{MR1031266}. 

The book of C.~Godsil and G.~Royle   \cite[Sec.15.13]{MR1829620} 
introduced a slightly different signed version of the Tutte polynomial,
$R(M;\alpha,\beta,x,y)$, following a result of K.~Murasugi \cite{MR0930077}. It is formulated in terms of signed matroids. Its application to the cycle matroid of a ribbon graph also can be obtained as a specialization $Z:=1$, $X:=x$, $Y:=y$ of our double weighted Bollob\'as-Riordan polynomial with weights of positive edges  $x_+:=\beta$, $y_+:=\alpha$, and weights of negative edges  $x_-:=\alpha$, $y_-:=\beta$.

\subsection{The dichromatic version of the Bollob\'as-Riordan polynomial.}

\def\bb{\mathbf b}
The dichromatic version of the Tutte polynomial $Z_G(a,b)$ (see for example 
\cite{MR0935433,MR1031266}), 
also known as the partition function of the Potts model in statistical mechanics, can be defined as a generating function of vertex colorings of a graph $G$ into $a$ colors:
$$Z_G(a,b)\ := \sum_{c\in Col(G)}\ 
(1+b)^{\mbox{\# edges colored not properly by $c$}}\ ,
$$
where $Col(G):=\{c:V(G)\to \{1,\dots,a\}\}$ is the set of colorings, and an edge $e$ is colored properly by $c$ if its endpoints are colored in different colors. Also, it can be expressed as a generating function of all spanning subgraphs:
$$Z_G(a,b) = \sum\limits_{F\subseteq E(\Ga)}a^{k(F)} b^{e(F)}.
$$
The Tutte polynomial is equivalent to the dichromatic
polynomial because of the identities
$$\begin{array}{rcl}
Z_G(a,b) &=& a^{k(G)}b^{r(G)} T(G;1+ab^{-1},1+b)\ ;\\
T(G;x,y) &=& (x-1)^{-k(G)}(y-1)^{-v(G)} Z_G((x-1)(y-1),y-1)\ .
\end{array}
$$

A multivariable version of dichromatic polynomial was introduced in \cite{MR0955462} 
and used in \cite{MR2187739}. 
The multivariable Tutte polynomial also appears as a very special case of Zaslavsky's colored Tutte polynomial of a matriod 
\cite{MR1080738}. 
It can be defined as 
$$Z_G(a,\bb):=\sum_{F\subseteq E(G)} a^{k(F)} \prod_{e\in F}b_e\ ,
$$
The sum runs over all spanning subgraphs of $G$, which we identify with subsets $F$ of $E(G)$; $\bb:=\{b_e\}$ is the set of variables (weights).

The multivariable dichromatic polynomial was generalized to ribbon graphs in  \cite{MR2368618} 
as
$$Z_G(a,\bb,c):=\sum_{F\subseteq E(G)} a^{k(F)} 
   \Bigl(\prod_{e\in F}b_e\Bigr) c^{bc(F)}.
\index{Bollob\'as--Riordan polynomial!dichromatic version}
$$
It is equivalent to the Bollob\'as-Riordan polynomial $BR_G(X,Y,Z)$ due to the relations
$$Z_G(a,\bb,c)= (ac)^{k(G)}BR_G(ac,c,c^{-1})\ ,
$$
$$BR_G(X,Y,Z)= \Bigl(\prod_{e\in E(G)}y_e\Bigr) (YZ)^{-v(G)}X^{-k(G)}
Z_G(XYZ^2,\{x_eYZ/y_e\},Z^{-1})\ ,
$$
where on the right-hand side of the first equation we use the Bollob\'as-Riordan polynomial with weights $x_e:=b_e$, $y_e:=1$.

\subsection{The arrow version of the Bollob\'as-Riordan polynomial.}

The arrow version of the Jones polynomial was introduced in \cite{DK} as far reaching generalization of the classical Jones polynomial for virtual links. The Thistlethwaite-type theorem for this polynomial was obtained in \cite{MR2994589}. 

\begin{definition}
An {\it arrow ribbon graph} is a ribbon graph together with
an {\it arrow structure} on it which is a (possibly empty) set of arrows tangent to the boundaries of the (vertex- and edge-) discs of the decomposition.
Two arrow graphs are {\it equivalent} if there is a homeomorphism between the corresponding surfaces respecting the decompositions and the orientations of the arrows.
The endpoints of segments along which the edges are attached to the vertices divide the boundaries of (vertex- and edge-) discs into arcs.
The arrows may slide along these arcs but may not slide over the end-points of the arcs or each other. Each arc may contain several arrows.
\end{definition}

So, if there are several arrows on an arc, then only their order on the arc is relevant, not their actual position on the arc. 

Here are some examples.
$$\mbox{(a)}\risS{-42}{ag-ex1}{}{40}{0}{0}\hspace{1.4cm}
  \mbox{(b)}\risS{-34}{ag-ex2}{}{40}{0}{0}\hspace{1.4cm}
  \mbox{(c)}\risS{-42}{ag-ex3}{}{75}{0}{40}
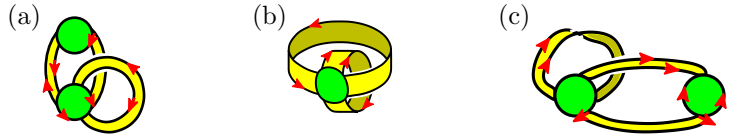
$$

\def\ap#1{\overrightarrow{BR}_{#1}}
In the presence of an arrow structure on a ribbon graph $G$ we can extend
the Bollob\'as-Riordan polynomial to the {\it arrow Bollob\'as-Riordan polynomial} 
\cite{MR2994589} 
adding new variables $K_c$ reflecting the arrow structure:
$$\ap{G}:=\!\!\sum_{F\subseteq E(G)} 
  (\prod_{e\in F}x_e)\ (\prod_{e\not\in F}y_e) 
  X^{r(G)-r(F)}Y^{n(F)}Z^{k(F)-bc(F)+n(F)}\!\!\!\!
   \prod_{f\in\partial(F)}\!\!\!\! K_{c(f)}.
\index{Bollob\'as--Riordan polynomial!arrow version}
$$
Here the independent variables $K_{c(f)}$ are assigned to each boundary component with 
$c(f)$ being equal to half of the number of arrows along the boundary component $f$ remaining after cancellations of all pairs of arrows which are next to each other and pointing in the same direction.
$$\risS{-25}{ar-can1}{}{50}{20}{30}\ \risS{-2}{toto}{}{25}{0}{0}\ 
      K_1\qquad
\risS{-25}{ar-can2}{}{50}{20}{30}\ \risS{-2}{toto}{}{25}{0}{0}\ 
      K_{1/2}\qquad
\risS{-25}{ar-can3}{}{50}{20}{30}\ \risS{-2}{toto}{}{25}{0}{0}\ 
      K_2
$$
We set $K_0=1$. One may note that whenever the number of arrows is odd on a boundary component, the associated variable is always $K_{1/2}$.
Thus the arrow Bollob\'as-Riordan polynomial $\ap{G}$ becomes a polynomial in infinitely many variables 
$x_e, y_e, X, Y, Z, K_{1/2}, K_1, K_2, \dots$.
Note that, for a concrete graph $G$ only finitely many $K$'s appear in $\ap{G}$.

\begin{example}
For the arrow graph $G$ shown on the leftmost column in the table, there are eight spanning subgraphs. Their parameters and the corresponding monomial in $K$'s are shown.

$$\begin{array}{c||c|c|c|c} \label{br-table}
\!\!\risS{8}{ag-ex3}{\put(23,30){$e_1$}\put(40,24){$e_2$}
                           \put(40,-5){$e_3$}}{56}{30}{5}
& \risS{8}{agBBB}{}{52}{0}{0}
& \risS{8}{agBBA}{}{52}{0}{0} 
& \risS{15}{agBAB}{}{52}{0}{0}  
& \risS{15}{agBAA}{}{52}{0}{0}\\ \hline
 (k,r,n,bc) & (1,1,1,2) & (1,1,0,1) & (1,1,0,1) & (2,0,0,2)\makebox(0,12){} \\
 \prod K_{c(f)} & K_1^2 & K_1 & K_1 & K_{1/2}^2\\ \hline\hline
& \!\risS{8}{agABB}{}{56}{40}{0}\!
& \!\risS{8}{agABA}{}{56}{0}{0}\!
& \!\risS{10}{agAAB}{}{56}{0}{0}\!
& \!\risS{10}{agAAA}{}{56}{0}{0}\!\!\!\!\!\\ \cline{2-5}
& (1,1,2,1) & (1,1,1,1) & (1,1,1,1) & (2,0,1,2)\makebox(0,12){}\\
& K_1 & K_1 & 1 &  K_{1/2}^2
\end{array}
$$
\end{example}
$$\begin{array}{ccl}
\ap{G} &=& y_1x_2x_3YK_1^2+y_1y_2x_3K_1+y_1x_2y_3K_1+
              y_1y_2y_3XK_{1/2}^2 \\[5pt]
&&\hspace*{-20pt}
 +x_1x_2x_3Y^2Z^2K_1+x_1y_2x_3YZK_1+x_1x_2y_3YZ+x_1y_2y_3XYZK_{1/2}^2\ .
\end{array}
$$

\section{The relative Tutte polynomial of planar graphs with distinguished edges.} 
\begin{definition}
A {\it relative 
graph} is a 
graph $G$ with a distinguished subset 
$H \subseteq E(G)$ of edges. The edges $H$ are called the 0-{\it edges} of $G$. Edges in $E(G)\backslash H$ will be referred to as {\it regular edges}. 
\end{definition}\index{relative graph}

We will use only plane relative graphs.

\begin{definition}
The {\it medial graph} $M(G)$, see for example \cite[Section 1.5.1]{MR3086663} or
 \cite[Section 7.7.1]{GrYe}, of a ribbon graph $G$ is a 4-valent graph embedded in the same surface as $G$ and constructed in a standard way. Consider the core graph $\Ga$ of $G$ as naturally embedded in $G$. Place a vertex at the center of each edge of $\Ga$. These will be the vertices of $M(G)$. The centers of edges of
$\Ga$ adjacent in the cyclic ordering around a vertex of $\Ga$ are connected by arcs following the boundary of $G$ and not intersecting $\Ga$. These will be the edges of $M(G)$. 
\end{definition}\index{medial graph}

Here is an example of the construction for a plane graph. 
$$\risS{-60}{med2}{}{250}{20}{70}
$$
The medial graph can be considered as an image of an immersion of several circles into the plane having only double points at the vertices of $M(G)$ as intersections (and self-intersections). The number of these circles is denoted by $\delta(G)$.

The last ingredient in definition of the relative Tutte polynomial is a function $\psi(G)$ of two variables $d$ and $w$ defined for plane ribbon graph $G$ as follows
$$\psi(G):=d^{\delta(G)-k(G)}w^{v(G)-k(G)}\ .
$$
This function is {\it block-invariant}, which means that its value on a one-point join of two graphs $G_1$ and $G_2$ does depend on the choice of identifying vertices of $G_1$ and $G_2$. 

\begin{definition}
Let $G$ be a relative plane graph with the distinguished set of 0-edges $H$.
We define the {\it relative Tutte polynomial}
\index{relative Tutte polynomial} by the equation
$$T_{G,H}(X,Y,\psi):=\sum_{F \subseteq G\setminus H} 
         (\prod_{e \in F}x_e) (\prod_{e \in \overline{F}}y_e) 
         X^{k(F\cup H)-k(G)}Y^{n(F)}\psi (H_F)\ ,
$$
here, abusing notation, we use $F$, $H$, $F \cup H$ etc. both for denoting the subsets of edges and the spanning subgraph of $G$ with those edges;  $\overline{F}:=G \setminus (F\cup H)$; and 
$H_F$ is the plane graph obtained from $F \cup H$ by contracting all edges of $F$. 
\end{definition}

\bigskip
{\bf Remarks.}\ 

{\bf 1.} The relative Tutte polynomial was introduced by Y.~Diao and G.~Hetyei in \cite{MR3111544}, 
who used the notion of {\it activities} to produce the most general form of it. The all subset formula we use was discovered by a group of undergraduate students (M.~Carnovale, Y.~Dong, J.~Jeffries) at the Ohio State University summer program ``Knots and Graphs" in 2009.
However, similar expressions may be traced back to L.~Traldi 
\cite{MR2047240} 
for the non-relative case, and to S.~Chaiken \cite{MR0982858} 
for the relative case of matroids. 

{\bf 2.} The function $\psi$ in \cite{MR3111544} can be obtained from ours by the substitution $w=1$. 

{\bf 3.} Another difference with \cite{MR3111544} is that we are using a doubly weighted version of the relative Tutte polynomial with weights $(x_e,y_e)$ of an edge $e\in G\setminus H$.

{\bf 4.} In the process of constructing the graph $H_F$ by contracting the edges of $F$ in $F \cup H$, we may come to a situation when we have to contract a loop. Then the contraction of a loop actually means its deletion. Since $G$ and $F \cup H$ are plane graphs, then the graph 
$H_F$ is also embedded in the plane.

{\bf 5.} While the medial graph of the planar graph $H_{F}$ depends on the embedding of $H_{F}$ in the plane, the number $\delta(H_F)$ does not (see \cite{MR3111544}). It depends only on the abstract graph $H_F$.
This was proved in the paper \cite{MR0587623}, 
which showed that $\delta(H_F)= \log_2|T(H_F;-1,-1)| + 1$. 

\subsection{From ribbon graphs to relative plane graphs.}

Here we describe a construction from \cite{MR3390227} 
relating ribbon graphs and planar relative graphs.

The manner in which we draw ribbon graphs suggests that we consider a projection $\pi: G\to \R^2$ with singularites of  two types. The first occurs when two edge-ribbons cross over each other. The second occurs when an edge ribbon twists over itself.  Away from these singularities the projection is one-to-one. In fact we use this idea on the medial graph $M(G)$.
$$
\risS{-30}{singular2}{\put(45,60){$R$}\put(25,40){$\pi$}
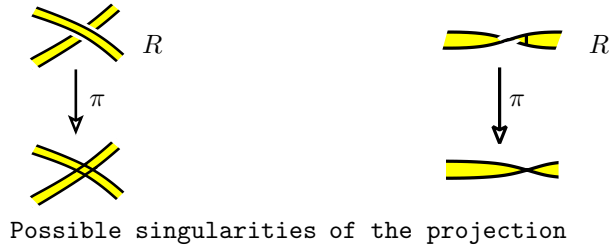
        \put(-5,-10){\tt Possible singularities of the projection}
         }{40}{0}{0}\hspace{4cm}
\risS{-20}{singular}{\put(60,50){$R$}\put(30,30){$\pi$}
        }{60}{50}{45}
$$
The image $B:=\pi(M(G))$ may then be considered as a regular 4-valent planar graph whose vertices are divided into two types. The vertices which are images of vertices of $M(G)$ will be called regular vertices, and the vertices that arise from the singularities of the projection will be called 0-vertices. 
The faces, the connected components of the complement of $B$ in the plane, can be colored in white and green in checkerboard manner.
Place a vertex on each green face and connect the vertices corresponding to the faces sharing a vertex of $B$ by an edge through the corresponding vertex of $B$. We obtain get a plane graph 
$C(B)$ whose medial graph is $B$, $M(C(B))=B$. This procedure is similar to taking the plane dual,
only now we place a vertex only on green faces, not on every face, and th edges of $C(B)$ go not across the edges of $B$ but across the vertices of $B$.
We consider $C(B)$ as a relative plane graph whose 0-edges are those which pass through the 
0-vertices of $B$.

The constructed relative plane graph $C(B)$ clearly depends on the projection $\pi$ and on the position of the vertices of the medial graph on the edge-ribbons.
However the invariants we will work with will not be affected by this ambiguity. The figure below shows two possibilities for $C(B)$ depending on the position of a vertex of the medial graph $M(G)$.

\begin{example}
$$\risS{-100}{rg1a}{\put(225,65){$B$} \put(285,65){$C(B)$}
                    \put(40,22){$\pi(G)$}\put(137,75){$M(G)$}
                    \put(53,135){$G$}\put(33,100){$\pi$}}{320}{35}{100}
$$
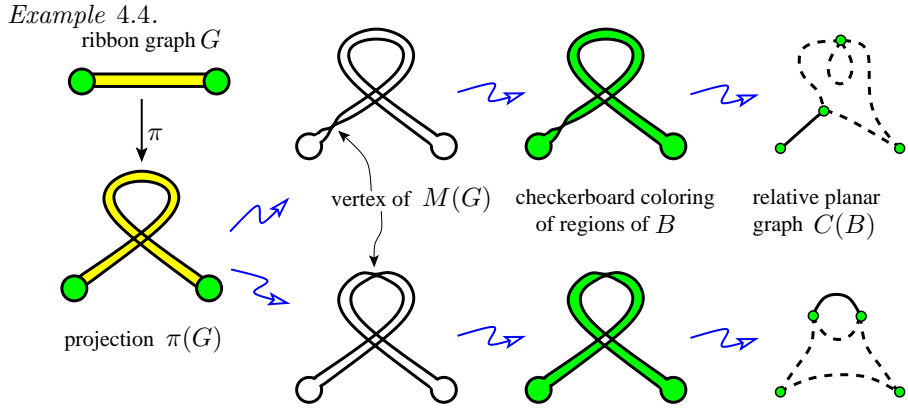
In this figure the 0-edges of $C(B)$ are drawn as dashed lines.
\end{example}

\subsection{From relative plane graphs to ribbon graphs.}\label{ss:rpg2rg}\ 

Conversely, from a relative plane graph $C$ we may construct a ribbon graph $G$. Consider the spanning subgraph $H$ of $C$ whose edges are the 0-edges of $C$. Construct its medial graph $M(H)$. Consider it as an immersion of a collection of $\d(H)$ circles with only transverse double points as singularities. 
We construct a ribbon graph $G$ with $\d(H)$ vertices corresponding to the circles of the medial graph $M(G)$. The edges of $G$ correspond to regular edges of $C$. They glued to the 
vertex-discs according to the intersection of a small neighborhood of a regular edge with the arcs of the medial graph. Namely, we can label small arcs on $M(G)$ at which the regular edges 
intersect $M(G)$ by arrows following counterclockwise orientation of the plane. Then we glue ribbons along the corresponding arcs of the vertex-discs.  


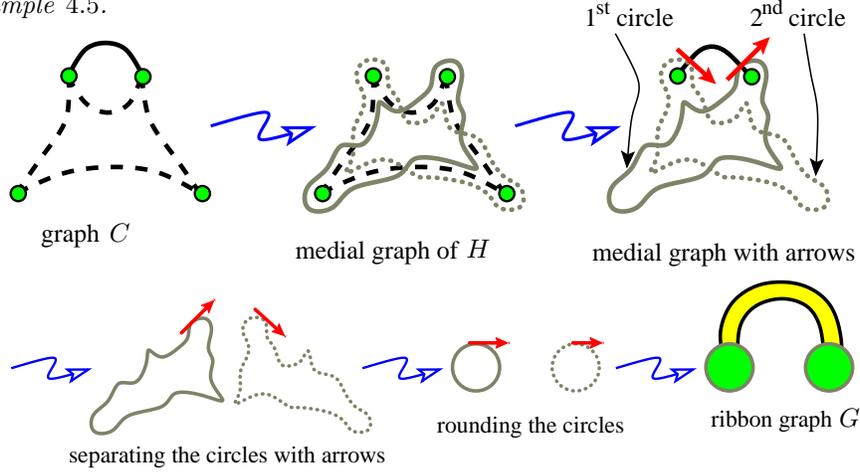
\begin{example}
$$\risS{-80}{pg-rg1}{\put(38,10){$C$} \put(173,4){$H$}
                     }{320}{0}{70}$$
$$\risS{-80}{pg-rg2}{\put(314,16){$G$}}{320}{0}{80}$$
\end{example} 

\subsection{Relation to the Bollob\'as--Riordan polynomial.} 

\begin{theorem}[\cite{MR3390227}]\label{th:BuCh}
{\it Suppose $G$ is a ribbon graph, and $C$ is a relative plane graph associated to a projection of $G$ with $H$ as a spanning subgraph corresponding to the zero-edges. Or, equivalently, assume $C$ is a relative plane graph and $G$ is the ribbon graph arising from $C$. 

Then under the substitution $w=\sqrt{\frac{X}{Y}}, d=\sqrt{XY}$,
$$X^\a Y^\b T_{C,H}(X,Y,\psi)=BR_G( X,Y,\frac{1}{\sqrt{XY}} )\ ,
$$
where $\a:=k(C)-k(G)-\b$ and $\b:=-\frac{1}{2}(v(G)-v(C))$.}
\end{theorem}

{\bf Remarks.}\  
{\bf 1.} It is an interesting consequence of this theorem that the specialization 
($w=\sqrt{\frac{X}{Y}}, d=\sqrt{XY}$) of the relative Tutte polynomial does not depend on the various  choices made in the construction of the relative plane graph. It is not difficult to describe a sequence of moves on relative plane graphs relating the graphs with different choices of the regular edges. It would be interesting to find such moves for different choices of the projection $\pi$ and, more generally, the moves preserving the relative Tutte polynomial.

{\bf 2.} The construction of $C$ from $G$ and back can be generalized to a wider class of projections $\pi$. We can require that only the restriction of $\pi$ to the boundary of $G$ be an immersion with only  ordinary double points as singularities. The theorem holds in this topologically more general situation. However, from the point of view of graph theory it is more natural to restrict ourselves to the class of projections above. 

\subsection{Dual relative plane graphs} \label{sss:dual}

Let $C$ be a relative plane graph with 0-edges $H$. The {\it dual} of $C$, denoted $C^*$ is formed by taking the dual of $C$ as a plane graph, and labeling the edges of $C^*$ which intersect 0-edges of $C$ as the 0-edges of $G^*$. Note that for relative plane graphs $(C^*)^*=C$, as with usual planar duality.\index{relative graph!dual}

The following theorem generalizes the classical relation for the Tutte polynomials of dual plane graphs, $T(C;x,y)=T(C^*;y,x)$, to relative plane graphs. 

\begin{theorem}[\cite{MR3390227}] 
Under the substitution $w=\sqrt{\frac{X}{Y}}$, $d=\sqrt{XY}$, we have
$$X^{a(C,H)}Y^{b(C)}T_{C,H}(X,Y,\psi)=
Y^{a(C^*,H^*)}X^{b(C^*)}T_{C^*,H^*}(Y,X,\psi)
$$
with the correspondence on the edge weights being $x_e=y_{e^*}$, $y_e=x_{e^*}$, where $e^*$ is the edge of $C^*$ that intersects $e$, and\qquad
$$a(C,H) = (|E(C \setminus H)|-v(C))/2+k(C)\ ,\qquad
b(C) = v(C)/2\ .
$$
\end{theorem}

\section{The Las Vergnas polynomial.} \label{s:LV}

\def\cC{\mathcal C}
M.~Las Vergnas \cite{MR0597150} came to his polynomials of graphs on surfaces through his general approach to matroid perspectives. It turns out that  for a (ribbon) graph on a surface with the edge set $E$, there is a geometric (Poincar\'e) dual graph $G^*$. For this graph there is a dual matroid to its cycle matroid,
$(\cC(G^*))^*$, which is the bond matroid of $G^*$. Then the natural bijection of the ground sets $(\cC(G^*))^*\to \cC(G)$ forms a matroid perspective according to \cite{MR0597150}. The Tutte polynomial of a perspective $M\to M'$ is a polynomial $T_{M\to M'}(x,y,z)$ in three variables defined by 
$$T_{M\to M'}\ := \sum_{F\subseteq M} (x-1)^{r_{M'}(E)-r_{M'}(F)} (y-1)^{n_M(F)}  
  z^{(r_M(E)-r_M(F))-(r_{M'}(E)-r_{M'}(F))},
$$
where in the rank and nullity functions, we use the subscript to indicate at which matroid the function fuction is considered in.

\bigskip
{\bf Properties} (\cite{MR0597150}). The usual Tutte polynomial of matroids $M$ and $M'$ can be recovered from the Tutte polynomial of matroid perspective in the following ways:\\
$\begin{array}{l}
T(M;x,y) = T_{M\to M}(x,y,z)\ ;\\
T(M;x,y) = T_{M\to M'}(x,y,x-1)\ ; \\
T(M';x,y) = (y-1)^{r(M)-r(M')} T_{M\to M'}(x,y,\frac1{y-1})\ . 
\end{array}$

\begin{definition}
The {\it Las Vergnas polynomial of $G$}, $LV_G(x,y,z)$,  is the Las Vergnas polynomial of the matroid perspective $(\cC(G^*))^*\to \cC(G)$.  \index{Las Vergnas polynomial}  
\end{definition}

If $G$ is a planar graph, then the Whitney planarity criteria claims that the matroid 
$(\cC(G^*))^*$ is isomorphic to $\cC(G)$. So the matroid perspective is the identity, and then the Las Vergnas polynomial is equal to the classical Tutte polynomial of $G$.

\begin{example}\label{ex14}
Let $G$ be a graph with one vertex and two loops embedded into a torus as shown. Then the dual graph $G^*$ is isomorphic to $G$.
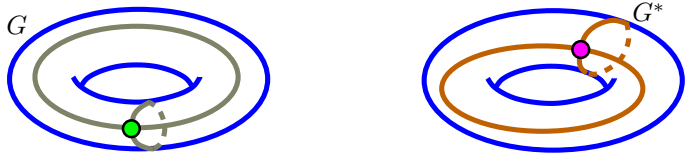
$$\risS{-20}{rg-torus}{\put(0,45){$G$}}{100}{0}{0}\hspace{2cm}
\risS{-20}{rg-torus-d}{\put(80,50){$G^*$}}{100}{35}{25}
$$
In this case the bond matroid $M=(\cC(G^*))^*$ has rank 2, and the cycle matroid $M'=\cC(G)$ has rank 0. 
For any subset $F$: $r_M(F)=|F|$, $n_M(F)=0$, and $r_{M'}(F)=0$. We have
$$LV_G(x,y,z)=z^2+2z+1\ .$$
\end{example}

\section{The Krushkal polynomial.} \label{s:Kp}

\def\S{\Sigma}
\def\k{\varkappa}
\def\im{\rm im}
The most general topological polynomial for graphs on surfaces was discovered by V.~Krushkal in his research in topological quantum field theories and the algebraic and combinatorial properties of models of statistical mechanics \cite{MR2769192}. 

Originally V.~Krushkal  \cite{MR2769192} defined his polynomial for orientable surfaces only. Here we follow the exposition of
Clark Butler \cite{But} who extended the definition to non-orientable surfaces.

\begin{definition}[\cite{But}]\label{d:kr-pol}
Let $G$ be a graph embedded in a surface $\S$. The embedding does not have to be cellular and the surface $\S$ does not have to be orientable.
Then the Krushkal polynomial is defined by
$$
K_{G,\S}(X, Y, A, B) := \sum_{F\subseteq G} X^{k(F)-k(G)}
Y^{k(\S\setminus F)-k(\S)} 
A^{s(F)/2} B^{s^\perp(F)/2}, \index{Krushkal polynomial}
$$
where the sum runs over all spanning subgraphs considered as ribbon graphs, and the parameters 
$s(F)$ and $s^\perp(F)$ are defined by
$$s(F):= 2k(F)-\chi(\widetilde{F}),\qquad
s^\perp(F):= 2k(\S\setminus F)-\chi(\widetilde{\S\setminus F})
$$
Here $\chi(\widetilde{F})$ (resp. $\chi(\widetilde{\S\setminus F})$)
stands for the Euler characteristic of the surface $\widetilde{F}$ 
(resp. $\widetilde{\S\setminus F}$) obtained by gluing a disc to each boundary component of the ribbon graph $F$ (resp. 
$\S\setminus F$).

As was outlined in Section \ref{s:BRpol}, for orientable surfaces the parameters $s$ and 
$s^\perp$ are equal to twice the genus of the corresponding surfaces, and for non-orientable surfaces they are equal to the number of M\"obius bands one has to glue into $k(F)$ and 
$\S\setminus F$ spheres to obtain a homeomerphic surfaces
\end{definition}

{\bf Remarks.}\  

{\bf 1.}  In the orientable case V.~Krushkal  \cite{MR2769192} used
a different expression for the exponent $\k$ of $Y$, as the dimension of the kernel of the map of the first homology groups
$$\k(F)=\dim(\ker(H_1(F;\R) \to H_1(\S;\R)))$$
induced by inclusion $F\hookrightarrow\S$. It was noted in \cite{MR2998838} that this dimension is equal to 
$\k(F)=k(\S\setminus F)-k(\S)$. C.~Butler observed in \cite{But} that considering homology groups with coefficients in $\Z_2:=\Z/2\Z$ leads to a topological interpretation of the exponents in the Krushkal polynomial for a general (not necessarily orientable) surface $\S$:
$$\begin{array}{rcl} 
k(\S\setminus F)-k(\S) &=&\dim(\ker(H_1(F;\Z_2) \to H_1(\S;\Z_2))),
\vspace{8pt}\\
s(F) &=&\dim H_1(\widetilde{F};\Z_2),\vspace{8pt}\\
s^\perp(F) &=&\dim H_1(\widetilde{\S\setminus F};\Z_2).
\end{array}
$$

{\bf 2.} In the orientable case V.~Krushkal \cite{MR2769192} indicated that the parameters $s(F)$ and $s^\perp(F)$ have another interpretation in terms of the symplectic bilinear form on the vector space $H_1(\S;\R)$ given by the intersection number.
For a given spanning subgraph $F$, let $V$ be its image in the homology group:
$$H_1(\S;\R)\supset V := V(F) := \im(H_1(F;\R) \to H_1(\S;\R))\ .
$$
For the subspace $V$ we can define its orthogonal complement $V^\perp$ in $H_1(\S;\R)$ with respect to the symplectic intersection form. Then
$$s(F)=\dim(V/(V\cap V^\perp))\ ,\qquad 
s^\perp(F)=\dim(V^\perp/(V\cap V^\perp))\ .
$$

\subsection{Relation to the previous polynomials.} 

Let $G$ be a ribbon graph, or equivalently a cellular embedding of the graph in a surface 
$\S:=\widetilde{G}$.

\begin{theorem}[\cite{MR2769192}]
The Tutte polynomial is a specialization of the Krushkal polynomial:
$$T(G;x,y)=(y-1)^{s(\S)/2} K_{G,\S}(x-1, y-1, y-1, (y-1)^{-1}).
$$
\end{theorem}

\begin{theorem}[\cite{MR2769192,But}]
The unweighted (with all edge weights equal to 1, $x_e=y_e=1$) Bollob\'as-Riordan polynomial is a specialization of the Krushkal polynomial:
\begin{equation}\label{eq:P2BR}
BR_G(X,Y,Z)=Y^{s(G)/2} K_{G,\S}(X,Y,YZ^2,Y^{-1}).
\end{equation}
\end{theorem}

\begin{theorem}[\cite{MR2998838,But}]
The Las Vergnas polynomial is a specialization of the Krushkal polynomial:
\begin{equation}\label{subst}
LV_G(x,y,z) = z^{s(G)/2} K_{G,\S}(x,y,z^{-1},z).
\end{equation}
\end{theorem}

\subsection{Properties.}
\begin{theorem}[\cite{MR2769192,But}]
$$K_{G,\S}=\left\{\begin{array}{ll}
K_{G/e,\S} + K_{G-e,\S} & 
         \mbox{if $e$ is ordinary, that is neither}\\
         &\mbox{a bridge nor a loop,} \\
(1+ X)\cdot K_{G/e} & \mbox{if $e$ is a bridge.} \\
(1+ Y)\cdot K_{G-e} & \mbox{if $e$ is a separable loop, the one whose removal}\\
         &\mbox{together with its vertex separates the surface $\S$.}
\end{array}\right.
$$
$$
K_{G_1\sqcup G_2,\S_1\sqcup \S_2} = K_{G_1,\S_1}\cdot K_{G_2,\S_2}, 
\mbox{where $\sqcup$ is the disjoint union.}
$$
\end{theorem}

\subsection{The contraction/deletion properties of these polynomials.}

The contraction/deletion properties for $BR_G$, for $K_{G,\S}$, and for $LV_{G,\S}$ are not quite the same. The problem arises when deletion of an edge of a ribbon graph changes its genus. The genus might decrease by 1 with the removal of an edge.
For example, if we delete a loop from the ribbon graph $G$ of Example \ref{ex14}, then the resulting graph with a single loop will have genus zero. So, while in the Bollob\'as-Riordan approach it is considered as a ribbon graph embedded into a sphere, in the Krushkal approach it is still embedded into the torus.
We cannot apply the substitution \eqref{eq:P2BR} to that graph
since its embedding on the torus is no longer cellular.
Thus the Krushkal polynomial does not satisfy the contraction/deletion property in the sense of the Bollob\'as and Riordan polynomial.
The Las Vergnas polynomial $LV_{G,\S}(x,y,z)$ does not satisfy either the contraction/deletion property in the sense of Bollob\'as and Riordan.
The contraction/deletion property for it discussed in \cite{EMM1}.

\begin{example}\label{ex3}
This is an example of a calculation of the three polynomials.
Here $G$ is a graph on the torus with two vertices and three edges $a$, $b$, and $c$ considered as a ribbon graph. Its dual $G^*$ has one vertex and three loops. We use the same symbols $a$, $b$, $c$ to denote the corresponding edges of $G^*$ as well.
$$
G=\risS{-22}{G-torus}{\put(52,33){$a$}\put(15,15){$a$}
      \put(46,8){$b$}\put(70,24){$b$}
      \put(28,13){$c$}\put(85,35){$c$}}{90}{30}{35}\ 
=\risS{-13}{cG-rg}{\put(37,32){$a$}\put(15,-7){$b$}
                        \put(33,-6){$c$}}{50}{0}{0}\qquad\qquad
G^*=\risS{-22}{Gd-torus}{\put(28,50){$a$}\put(72,16){$b$}
                        \put(7,26){$c$}}{90}{35}{35}\ 
$$
The matroid $M'=\cC(G)$ is of rank 1, and for any nonempty subset $F$, $r_{M'}(F)=1$. The cycle matroid $\cC(G^*)$ of the dual graph is of rank zero because $G^*$ has only loops. So its dual $M=(\cC(G^*))^*$ has rank 3, all subsets $F$ are independent and $r_M(F)=|F|$. The next table shows the values of various parameters and the contributions of all eight subsets $F\subseteq\{a,b,c\}$ to the three polynomials.
$$\begin{array}{cr||@{\ \makebox(0,10){}}c|c|c|c|c|c|c|c}
 &F& \emptyset&\{a\}&\{b\}&\{a,b\} &\{c\}&\{a,c\}&\{b,c\}&\{a,b,c\}\!\!\!  
        \\[2pt]
\hline\hline
 &k(F)& 2&1&1&1 &1&1&1&1 \\ \hhline{~---------}
 &\k(F)& 0&0&0&0 &0&0&0&0 \\ \hhline{~---------}
 &s(F)& 0&0&0&0 &0&0&0&2 \\ \hhline{~---------}
 &s^\perp(F)& 2&2&2&0 &2&0&0&0 \\ \hhline{~---------}
\makebox(20,0){\raisebox{60pt}{\rotatebox{90}{Krushkal}}}
 &K_{G,\S}& XB&B&B&1 &B&1&1&A \\ \hline\hline

 &r_M(F)& 0&1&1&2 &1&2&2&3 \\ \hhline{~---------}
 &r_{M'}(F)& 0&1&1&1 &1&1&1&1 \\ \hhline{~---------}
 &n_M(F)& 0&0&0&0 &0&0&0&0 \\ \hhline{~---------}
\makebox(20,0){\raisebox{50pt}{\rotatebox{90}{Las Vergnas}}}
 &LV_G& \!(x-1)z^2&z^2&z^2&z &z^2&z&z&1 \\ \hline\hline

 &k(F)& 2&1&1&1 &1&1&1&1 \\ \hhline{~---------}
 &n(F)& 0&0&0&1 &0&1&1&2 \\ \hhline{~---------}
 &bc(F)& 2&1&1&2 &1&2&2&1 \\ \hhline{~---------}
\makebox(20,0){\raisebox{50pt}{\rotatebox{90}{
              $\genfrac{}{}{0pt}{0}{\mbox{Bollob\'as}}{\mbox{Riordan}}$}}}
 &BR_G& X&1&1&Y &1&Y&Y&Y^2Z^2 \\ 
\end{array}
$$
Thus
$$K_{G,\S}=3+3B+XB +A,\qquad
LV_G=3z+3z^2+(x-1)z^2+1,
$$
$$
BR_G=3Y+3+X+Y^2Z^2.
$$
One can readily confirm the relations \eqref{subst} and \eqref{eq:P2BR} from here.

Now if we contract the edge $c$, the graph $G/c$ will still be cellularly  embedded in the same torus $\S$.
The right four columns of the table above give the following polynomials:
$$K_{G/c,\S}=B+2+A,\qquad
LV_{G/c}=z^2+2z+1,\qquad
BR_{G/c}=1+2Y+Y^2Z^2.
$$
Meanwhile if we delete the edge $c$, then 
$$K_{G-c,\S}=XB+2B+1\ .
$$
But the graph $G-c$ is not cellularly embedded into the torus $\S$ any more. Thus the Las Vergnas and the Bollob\'as-Riordan polynomials are not defined for it. The ribbon graph $G-c$, after gluing discs to its two boundary components, results in the sphere $\widetilde{G-c}=S^2$. Thus the graph $G-c$ embeds cellularly into the sphere $S^2$. For this embedding we have
$$K_{G-c,S^2}=X+2+Y,\qquad
LV_{G-c}=(x-1)+2+(y-1),\qquad
BR_{G-c}=X+2+Y.
$$
Therefore 
$$K_{G,\S}=K_{G-c,\S}+K_{G/c,\S}\qquad\mbox{and}\qquad
BR_{G}=BR_{G-c}+BR_{G/c}\ ,
$$
but
$$K_{G,\S}\not=K_{G-c,S^2}+K_{G/c,\S}\qquad\mbox{and}\qquad
LV_G\not=LV_{G-c}+LV_{G/c}\ .
$$
\end{example}

\subsection{Duality for the Krushkal polynomial.} \label{ss:dualK}

The classical relation, $T(G;x,y)=T(G^*;y,x)$, for the Tutte polynomials of dual plane graphs can be generalized to a relation for the Krushkal polynomial of arbitrary ribbon graphs. 

\begin{theorem}[\cite{MR2769192,But}]
Let $G$ be a ribbon graph naturally embedded into the surface $\S=\widetilde{G}$, and $G^*\subset\S$ be its Poincar\'e (geometric) dual ribbon graph.Then
$$K_{G,\S}(X, Y, A, B)=K_{G^*,\S}(Y, X, B, A)\ .
$$\index{Krushkal polynomial!duality}
\end{theorem}

\section{Quasi-trees.} \label{s:qt}

\begin{definition}[\cite{MR2854567}]
A {\it quasi-tree} is a ribbon graph with one boundary component.
\index{quasi-tree}
\end{definition}
 
Quasi-trees, as well as a notion of activities relative to a spanning quasi-tree, were 
introduced and used for a {\it quasi-tree expansion} of the Bollobas-Riordan polynomial in 
\cite{MR2854567} 
Their results were generalized to non-orientable ribbon graphs in \cite{MR2780852} 
(see also an unpublished manuscript \cite{Dew}). 

We denote the set of all spanning quasi-trees of a given ribbon graph $G$ by ${\mathcal Q}_G$. 

If $G$ is a plane graph of genus zero, then a spanning quasi-tree is a spanning tree. But for non-plane graphs the set of spanning quasi-trees is bigger than the set of spanning trees. In Example \ref{ex14}, the whole graph consisting of a vertex and two loops on a torus is a quasi-tree but, of course, not a tree. The table in Example \ref{ex3} shows that the ribbon graph $G$ has four spanning quasi-trees $\{a\}$, $\{b\}$, $\{c\}$, and $\{a,b,c\}$.

\subsection{Quasi-tree activities.} \label{ss:qt-act}

Let $\prec$ be a total order on edges $E(G)$ of a ribbon graph $G$, and $Q$ be a spanning quasi-tree of $G$. Tracing the boundary component of $Q$ we will get a round trip passing the boundary arcs of each edge-ribbon twice, on opposite sides of a ribbon separating it from faces for edges in $Q$ and on opposite sides of a ribbon separating it from vertices for edges not in $Q$. This can be encoded in a chord diagram $C_G(Q)$ consisting of a circle corresponding to the boundary of $Q$ with pairs of arcs corresponding to the same edge-ribbon connected by  chords.  Thus the set of chords inherits the total order $\prec$.

\begin{definition}[\cite{MR2854567}]
An edge is called {\it live} if the corresponding chord is smaller than any chord intersecting it relative to the order $\prec$. 
Otherwise it is called {\it dead}. \index{edge!live/dead}
\index{activities! w.r.t. quasi-tree}
\index{quasi-tree!activities}
\end{definition}

\begin{remark}
For a plane graph $G$ a spanning quasi-tree is a spanning tree and the notion of {\it live/dead} consides with the classical Tutte's notion of
{\it active/inactive}. But for the higher genus the notions are different, so we follow the terminology {\it live/dead} of \cite{MR2854567}, 
where a reader can find further discussion and examples.
Also, this notion is different from the {\it embedding-activities} introduced by O.~Bernardi in \cite{MR2428901}. 
\end{remark}

\begin{example}\label{ex4}
In Example \ref{ex14} of a single vertex and two loops there are two spanning quasi-trees consisting of a vertex without edges and the whole graph. Both of them have a chord diagram depicted on the left in the picture below.
In Example \ref{ex3} there are four spanning quasi-trees. All of them have the same chord diagram, depicted on the right in the picture below.
$$\risS{-10}{cd-2}{}{40}{30}{0}\qquad\qquad
\risS{-10}{cd-3}{}{40}{0}{0}
$$
In either case there is only one live edge, the smallest in the order $\prec$. It could be either {\it internal} or {\it external} depending on whether it belongs or does not belong to the spanning quasi-tree.
\end{example}

\begin{example}\label{ex5}
Here we consider a more complicated example of a non-orientable ribbon graph $G$ following 
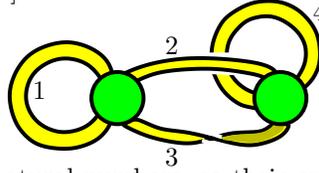
\cite{But}.
$$\risS{-20}{rg-ex4}{\put(10,20){$1$}\put(60,37){$2$}
               \put(60,-5){$3$}\put(116,50){$4$}}{120}{30}{20}
$$
The edges are labeled by natural numbers, so their order is the obvious one: 
$1\prec 2\prec 3\prec 4$.  There are four spanning quasi-trees 
$Q_{\{2\}}$, $Q_{\{3\}}$, $Q_{\{2,3\}}$, and $Q_{\{2,3,4\}}$, where the subscript indicates the set of edges included in the quasi-tree. Their chord diagrams look like this:
$$C_G(Q_{\{2\}})=
\risS{-17}{cd-4-q2}{\put(18,28){$\scriptstyle 1$}
              \put(28,22){$\scriptstyle 2$}
              \put(11,29){$\scriptstyle 3$}
              \put(26,12){$\scriptstyle 4$}}{40}{20}{30},\qquad 
C_G(Q_{\{3\}})=\risS{-17}{cd-4-q3}{\put(18,28){$\scriptstyle 1$}
              \put(32,21){$\scriptstyle 2$}
              \put(11,30){$\scriptstyle 3$}
              \put(20,4){$\scriptstyle 4$}}{40}{0}{0},
$$
$$C_G(Q_{\{2,3\}})=\risS{-17}{cd-4-q23}{\put(17,28){$\scriptstyle 1$}
              \put(28,20){$\scriptstyle 2$}
              \put(11,27){$\scriptstyle 3$}
              \put(4,12){$\scriptstyle 4$}}{40}{0}{20},\qquad 
C_G(Q_{\{2,3,4\}})=\risS{-17}{cd-4-q234}{\put(17,28){$\scriptstyle 1$}
              \put(28,19){$\scriptstyle 2$}
              \put(11,27){$\scriptstyle 3$}
              \put(14,8){$\scriptstyle 4$}}{40}{0}{0}.
$$
Thus $1$ is always an externally live edge, $2$ is also live, and $3$ is internally live for $Q_{\{2,3,4\}}$. We depict the chords corresponding to dead edges by dashed lines. Also some chords on the pictures are marked by a cross. These correspond to {\it non-orientable} edges according to the next definition. 
\end{example}

\begin{definition}[\cite{MR2780852}] 
An edge of a non-orientable ribbon graph $G$ can be classified as {\it orientable or 
non-orientable} relative to a spanning quasi-tree $Q$ for $G$.
Choose an orientation on the boundary circle $\partial Q$ of $Q$ and on the boundary of each edge-ribbon considered as a topological disc.
Then, tracing $\partial Q$ in the direction of the orientation we meet exactly two arcs of each-ribbon edge. If the orientations of these two arcs are coherent, either both along the tracing, or both against the tracing, then we call the edge {\it orientable} relative to 
$Q$. Otherwise we call it {\it non-orientable}.
\index{edge!orientable/non-orientable}

It is easy to see that this property does not depend on the choice of the orientations. Essentially the chord diagram $C_G(Q)$ depicts a partial dual ribbon graph $G^{E(Q)}$ in the sense of \cite{MR2507944}. All edges of $G^{E(Q)}$ are loops. For loops we can define orientability in the usual topological sense which coincides with our definition above.
\end{definition}

For the graph $G$ in Example \ref{ex5} consider a spanning quasi-tree
 and choose the orientations as follows.
$$\risS{-20}{rg-ex4}{\put(10,20){$1$}\put(60,37){$2$}
               \put(60,-5){$3$}\put(116,50){$4$}}{120}{40}{20}
\ \risS{-2}{toto}{}{25}{0}{0}\ 
\risS{-20}{rg-ex4or}{\put(14,20){$1$}\put(60,40){$2$}
               \put(60,-4){$3$}\put(128,50){$4$}}{130}{0}{0}
$$
We observe that the two relevant arcs of edge 1 are oriented coherently to
the direction of $\partial Q_{\{2,3\}}$. For the other edges this is not true.
So 1 is an orientable edge relative to $Q_{\{2,3\}}$, and the other three
edges are not.

Thus each edge has three bits of additional information relative to a spanning tree $Q$: internal/external, live/dead, and orientable/non-orientable. This information is used in the quasi-tree expansion of the Krushkal polynomial in next subsection.

\subsection{Quasi-tree expansion of the Krushkal polynomial.} \label{ss:qt-Kp}

\begin{definition}\label{d:FGamma}
With each spanning quasi-tree $Q$ of a ribbon graph $G$ we associate 
\begin{itemize}
\item[$\bullet$] a spanning ribbon subgraph $F(Q)$ by deleting the (internally) live orientable edges of $Q$;
\item[$\bullet$] an abstract (not embedded) graph $\Gamma(Q)$ whose vertices are the connected components of $F(Q)$ and whose edges are the internally live orientable edges of $Q$.
\end{itemize}
\end{definition}

We will also need a dual construction. For a ribbon graph $G$, regarded as a graph cellularly embedded into the surface $\S=\widetilde{G}$, we consider the Poincar\'e dual graph ribbon graph $G^*$. A spanning subgraph $F$ for $G$ determines a spanning subgraph $F^\star$ containing all edges of $G^*$ which do not intersect edges of $F$. Note that $F^\star$ is not a dual ribbon graph for $F$, its edge set is rather the complement of the set of edges of $G^*$ corresponding to the edges of $F$. The asterisk applied to an entire graph means the dual graph, but the star applied to a subgraph means the spanning subgraph given by this construction.

The spanning subgraphs $F$ and $F^\star$ have common boundary and their gluing along this common boundary gives the whole surface $\S$.
In particular, for a spanning quasi-tree $Q$ for $G$, the subgraph $Q^\star$ is a quasi-tree for $G^*$. Moreover, these quasi-trees have the same chord diagrams, $C_G(Q)=C_{G^*}(Q^\star)$. Also the natural bijection of edges of $G$ and $G^*$ leads to the total order $\prec^*$ on edges of 
$G^*$ induced by $\prec$. Consequently the property of an edge of being live/dead relative to 
$Q$ is mappeded by the bijection to the same property relative to $Q^\star$. Also one may check that the property of an edge of being orientable/non-orientable relative to $Q$ is preserved by the bijection to the same property relative to $Q^\star$. But its property of being internal/external is changed to the opposite.

\begin{definition}\label{d:F*Gamma*}
Now we can apply definition \ref{d:FGamma} to the quasi-tree of $G^*$ just constructed.
\begin{itemize}
\item[$\bullet$] We obtain a spanning ribbon subgraph $F(Q^\star)$ by deleting the internally live orientable edges of $Q^\star$. They correspond to externally live orientable edges of $Q$.
\item[$\bullet$] We obtain an abstract graph $\Gamma(Q^\star)$ consists of vertices that are the connected components of $F(Q^\star)$ and edges that are the internally live orientable edges of 
$Q^\star$.
\end{itemize}
\end{definition}

\begin{theorem}[\cite{But}]\label{th:But}
For a ribbon graph $G$, the Krushkal polynomial has the following expansion over the set of quasi-trees. 
$$K_G(X,Y,A,B)=\sum_{Q\in{\mathcal Q}_G}
A^{s(F(Q))/2}T_Q \cdot B^{s(F(Q^\star))/2}T_{Q^\star}\ ,
$$
where $T_Q=T(\Gamma(Q);X+1,A+1)$ and $T_{Q^\star}=T(\Gamma(Q^\star);Y+1,B+1)$ stand for the classical Tutte polynomial of the abstract graphs
$\Gamma(Q)$ and $\Gamma(Q^\star)$, and the parameter $s$ was defined in \ref{d:kr-pol}.
\index{Krushkal polynomial!quasi-tree expansion}
\end{theorem}

\begin{example}\label{ex6}
For the graph $G$ on the torus with a single vertex and two loops, parallel and meridian, in Example \ref{ex14} we have two spanning quasi-trees consisting of a vertex only and the whole graph, as in Example \ref{ex4}.
For $Q$ being a vertex, both $F(Q)$ and $\Gamma(Q)$ consist of one vertex. So $s(F(Q))=0$ and 
$T_Q=1$. The dual quasi-tree $Q^\star=D^*$ is the whole graph in this case, and one of its loops is internally active. 
So $F(Q^\star)$ consists of one vertex with a loop attached to it. 
The abstract graph $\Gamma(Q^\star)$ is the same. Thus $s(F(Q^\star))=0$ and 
$T_{Q^\star}=T(\Gamma(Q^\star);Y+1,B+1)=B+1$. Therefore the contribution of the single vertex spanning quasi-tree is $B+1$.
For the second quasi-tree $Q$, the whole graph, the situation is completely symmetrical. Its dual quasi-tree $Q^\star$ consists of a single vertex, and its contribution is $A+1$. So
$$K_G(X,Y,A,B)=A+2+B,
$$
as was mentioned in Example \ref{ex3}.
\end{example}

\begin{example}\label{ex7}
For the graph $G$ in Examples \ref{ex3} and \ref{ex4} with edges labeled $a$, $b$, and $c$ there are four spanning quasi-trees $Q_{\{a\}}$, $Q_{\{b\}}$, $Q_{\{c\}}$, and $Q_{\{a,b,c\}}$.
For each of them $a$ is a live edge and $b$ and $c$ are dead relative to the total order $a\prec b\prec c$. The next tables show the ribbon graphs
$F(Q)$, $F(Q^\star)$ and abstract graphs $\Gamma(Q)$, $\Gamma(Q^\star)$ and their contributions to the Krushkal polynomial for each of these spanning quasi-trees.
$$\begin{array}{@{\rb{-6pt}{\makebox(0,15){}}}c||c|c|c|c} \hline\hline
Q&Q_{\{a\}}&Q_{\{b\}}&Q_{\{c\}}&Q_{\{a,b,c\}} \\ \hline\hline
F(Q)&\risS{-5}{F_Q_a}{}{50}{15}{25}&
     \risS{-13}{F_Q_b}{\put(15,-7){$b$}}{50}{0}{0}&
     \risS{-13}{F_Q_c}{\put(33,-6){$c$}}{50}{0}{0}&
     \risS{-13}{F_Q_abc}{\put(15,-7){$b$}\put(33,-6){$c$}}{50}{0}{0} 
     \\ \hline
\makebox(0,15){} A^{s(F(Q))/2}&1&1&1&1\\ \hline
\Gamma(Q)&\risS{0}{Ga_Q_a}{}{30}{20}{15}&
     \risS{0}{Ga_Q_b}{}{6}{0}{0}&\risS{0}{Ga_Q_b}{}{6}{0}{0}&
     \risS{-10}{Ga_Q_abc}{}{20}{0}{0} \\ \hline
T_Q&X+1&1&1&A+1 \\ \hline \hline
&&&&\\ \hline\hline
Q^\star&Q^\star_{\{a\}}&Q^\star_{\{b\}}&Q^\star_{\{c\}}&Q^\star_{\{a,b,c\}} \\ \hline\hline
F(Q^\star)&\risS{-12}{F_Qst_a}{\put(15,-7){$b$}\put(33,-6){$c$}}{50}{25}{25}&
     \risS{-12}{F_Qst_b}{\put(15,-5){$c$}}{30}{0}{0}&
     \risS{-12}{F_Qst_c}{\put(20,-6){$b$}}{30}{0}{0}&
     \risS{-2}{F_Qst_abc}{}{14}{0}{0} 
     \\ \hline
\makebox(0,15){} B^{s(F(Q^\star))/2}&B&1&1&1\\ \hline
\Gamma(Q^\star)&\risS{0}{Ga_Q_b}{}{6}{20}{15}&
     \risS{-10}{Ga_Q_abc}{}{20}{0}{0}&\risS{-10}{Ga_Q_abc}{}{20}{0}{0}&
     \risS{0}{Ga_Q_b}{}{6}{0}{0} \\ \hline
T_{Q^\star}&1&B+1&B+1&1 \\ \hline \hline
\end{array}
$$
So the Krushkal polynomial is equal to
$$K_G=(X+1)B+(B+1)+(B+1)+(A+1)=XB+A+3B+3,
$$
which coincides with the direct calculation in Example \ref{ex3}
\end{example}

\begin{example}\label{ex8}
Now we consider the example \ref{ex5} of a non-orientable ribbon graph
$G$ and its dual $G^*$ from \cite{But}.
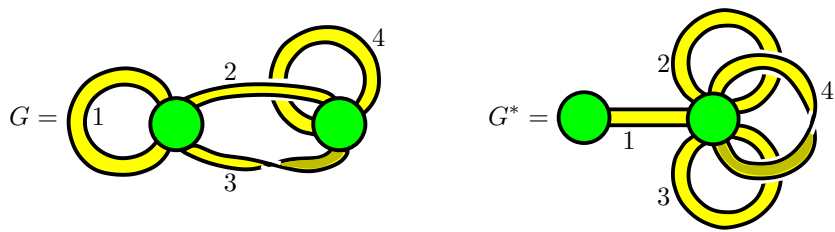
$$G=\risS{-20}{rg-ex4}{\put(10,20){$1$}\put(60,37){$2$}
               \put(60,-5){$3$}\put(116,50){$4$}}{120}{30}{35}
\qquad\qquad
G^*=\risS{-40}{rgd-ex4}{\put(25,31){$1$}\put(38,60){$2$}
               \put(38,10){$3$}\put(100,50){$4$}}{100}{0}{0}
$$
It has four spanning quasi-trees $Q_{\{2\}}$, $Q_{\{3\}}$, $Q_{\{2,3\}}$, and $Q_{\{2,3,4\}}$. Again we indicate the ribbon graphs
$F(Q)$, $F(Q^\star)$ and abstract graphs $\Gamma(Q)$, $\Gamma(Q^\star)$ and their contributions to the Krushkal polynomial.
$$\begin{array}{@{\rb{-6pt}{\makebox(0,15){}}}c||c|c|c|c} \hline\hline
Q&Q_{\{2\}}&Q_{\{3\}}&Q_{\{2,3\}}&Q_{\{2,3,4\}} \\ \hline\hline
F(Q)&\risS{-5}{F_Q_a}{}{50}{15}{25}&
     \risS{-10}{F_Q_3}{\put(15,-7){$3$}}{50}{0}{0}&
     \risS{-16}{F_Q_23}{\put(30,22){$2$}\put(15,-6){$3$}}{50}{0}{0}&
     \risS{-20}{F_Q_234}{\put(55,20){$4$}}{55}{0}{0} 
     \\ \hline
\makebox(0,15){} A^{s(F(Q))/2}&1&1&A^{1/2}&1\\ \hline
\Gamma(Q)&\risS{0}{Ga_Q_a}{}{30}{20}{15}&
     \risS{0}{Ga_Q_b}{}{6}{0}{0}&\risS{0}{Ga_Q_b}{}{6}{0}{0}&
     \risS{-5}{Ga_Q_234}{}{35}{0}{0} \\ \hline
T_Q&X+1&1&1&X+A+2 \\ \hline \hline
&&&&\\ \hline\hline
Q^\star&Q^\star_{\{2\}}&Q^\star_{\{3\}}&Q^\star_{\{2,3\}}&Q^\star_{\{2,3,4\}} \\ \hline\hline
F(Q^\star)&
     \risS{-18}{F_Qst_2}{\put(18,6){$3$}\put(57,26){$4$}}{57}{25}{25}\ &
     \risS{-18}{F_Qst_3}{\put(22,32){$2$}\put(53,26){$4$}}{57}{0}{0}&
     \risS{-12}{F_Qst_23}{\put(53,26){$4$}}{57}{0}{0}&
     \risS{-5}{F_Q_a}{}{50}{0}{0} 
     \\ \hline
\makebox(0,15){} B^{s(F(Q^\star))/2}&B&B&B^{1/2}&1\\ \hline
\Gamma(Q^\star)&\risS{0}{Ga_Q_a}{}{30}{20}{15}&
     \risS{0}{Ga_Q_a}{}{33}{0}{0}&\risS{0}{Ga_Q_a}{}{30}{20}{15}&
     \risS{0}{Ga_Q_a}{}{30}{20}{15} \\ \hline
T_{Q^\star}&Y+1&Y+1&Y+1&Y+1 \\ \hline \hline
\end{array}
$$
So the Krushkal polynomial is equal to
$$K_G=(X+1)(Y+1)B+(Y+1)B+(Y+1)A^{1/2}B^{1/2}+(X+A+2)(Y+1).
$$
One may check this answer by a direct calculation according to the definition \ref{d:kr-pol}.
\end{example}

\bigskip
Clark Butler's Theorem \ref{th:But} specializes to the quasi-tree expansion of the Bollob\'as-Riordan polynomial from \cite{MR2854567} 
in the orientable case and from \cite{MR2780852} 
in the non-orientable case.
Also it specializes to a quasi-tree expansion of the Las Vergnas polynomial, see \cite{But}. It would be very interesting to find such an expansion for the Tutte polynomial of matroid perspectives entirely in the matroid setting.


\end{document}